# Global weak solutions for a two-species Riesz cross-attraction system with quadratic diffusion

**Zhenxing Liu**[1], **Shen Bian**[1,*]

[1] Beijing University of Chemical Technology, Beijing 100029, China

* Corresponding author: 2024201118@buct.edu.cn

**Abstract.** We establish global existence for a two-species cross-attraction system on $\mathbb{R}^3$, with quadratic diffusion and the common interaction operator $(-\Delta)^{-s}$, where $1/2 < s < 1$. Nonnegative initial data in $L^1 \cap L^\infty$ with finite second moments require no smallness condition on the masses. A joint energy and second-moment estimate and a coupled $L^p$ iteration give uniform a priori estimates for the regularized solutions. Compactness then yields mass-conserving global weak solutions, bounded on every finite time interval and strongly continuous in $L^1 \cap L^2$. The energy inequality holds at every pair of time endpoints. Uniqueness and bounds uniform for all positive times are not asserted.

## 1. Introduction

We study the following Cauchy problem for two mutually attracting densities on $\mathbb{R}^3 \times (0,\infty)$:

$$\begin{aligned} &\partial_t u = \Delta(u^2) - \kappa_1 \nabla \cdot \big(u \nabla (K_s * w)\big), \\ &\partial_t w = \Delta(w^2) - \kappa_2 \nabla \cdot \big(w \nabla (K_s * u)\big), \\ &u(x,0) = u_0(x), \quad w(x,0) = w_0(x). \end{aligned} \tag{1}$$

Here $\kappa_1, \kappa_2 > 0$ and $K_s(x) = c_{3,s} |x|^{-(3-2s)}$ for $x \neq 0$, where $c_{3,s} = \Gamma((3-2s)/2)/(2^{2s}\pi^{3/2}\Gamma(s))$ and $1/2 < s < 1$. The convolution operator is the Riesz integral of order $2s$; the convolution is over $\mathbb{R}^3$. The diffusion is local and of porous-medium type.

Keller–Segel theory motivates aggregation–diffusion models.[1,2] For smooth kernels and linear diffusion, Giunta et al. prove local existence for multispecies systems and global existence in one dimension.[3] Further studies treat multispecies singular kernels,[4] degenerate multipopulation systems,[5] and irregular kernels.[6] The closest comparisons are Newtonian two-species cross-attraction[7] and scalar Riesz interactions.[8] Related studies concern nonlocal-to-local limits for quadratic multispecies cross-diffusion[9] and equilibria associated with nonlinear Riesz potentials.[10]

The gradient $|\nabla K_s(x)| = (3-2s)c_{3,s}|x|^{-(4-2s)}$ is more singular near zero than the Newtonian order $|x|^{-2}$. In Lemma 2.3, Step 2 (p. 10), Carrillo and Lin use $-\Delta v_\varepsilon = w_\varepsilon$ to obtain $\int_{\mathbb{R}^d} u_\varepsilon^p w_\varepsilon \, dx$ in (2.16).[7] Here, for smooth data, $-\Delta(K_s * w) = (-\Delta)^{1-s} w$, so their local-product estimate does not transfer directly. Lemma 3.3 instead uses HLS estimates of order $2s-1$ and joint absorption by both diffusion dissipations. Coupled testing is already used by Carrillo and Lin[7]; the new estimates handle the fractional interaction.

We use fractional integration,[11–13] quasilinear existence,[14] parabolic regularity,[15,16] and compactness.[17] Energy and second moments provide quadratic bounds, followed by the coupled iteration. Constants concern fixed $s \in (1/2,1)$, without uniformity as $s \downarrow 1/2$.

## 2. Main result and energy structure

Assume $u_0, w_0 \ge 0$, $u_0, w_0 \in L^1(\mathbb{R}^3) \cap L^\infty(\mathbb{R}^3)$, with positive masses $M_u = \int_{\mathbb{R}^3} u_0(x)\,dx$ and $M_w = \int_{\mathbb{R}^3} w_0(x)\,dx$, and $\int_{\mathbb{R}^3} |x|^2 (u_0 + w_0)(x)\,dx < \infty$. Write $Q_T = \mathbb{R}^3 \times (0,T)$ and $\|\cdot\|_p = \|\cdot\|_{L^p(\mathbb{R}^3)}$. Here $B_L = \{x \in \mathbb{R}^3 : |x| < L\}$. Unspecified spatial function spaces are over $\mathbb{R}^3$; Sobolev spaces are nonhomogeneous.

**Definition 2.1.** A nonnegative pair is a weak solution of (1) if the terms below are locally integrable and, for every $\phi \in C_c^\infty(\mathbb{R}^3 \times [0,T))$, both identities hold:

$$\begin{aligned} &\int_0^T\int_{\mathbb{R}^3} \left(u\phi_t + u^2\Delta\phi + \kappa_1 u \nabla(K_s * w)\cdot\nabla\phi\right)dxdt \\ &\qquad + \int_{\mathbb{R}^3} u_0(x)\phi(x,0)\,dx = 0, \\ &\int_0^T\int_{\mathbb{R}^3} \left(w\phi_t + w^2\Delta\phi + \kappa_2 w \nabla(K_s * u)\cdot\nabla\phi\right)dxdt \\ &\qquad + \int_{\mathbb{R}^3} w_0(x)\phi(x,0)\,dx = 0. \end{aligned} \tag{2}$$

Define the energy and weighted fluxes by

$$E(u,w) = \frac{1}{\kappa_1}\int_{\mathbb{R}^3} u^2\,dx + \frac{1}{\kappa_2}\int_{\mathbb{R}^3} w^2\,dx - \int_{\mathbb{R}^3} u(K_s * w)\,dx, \tag{3}$$

$$\begin{aligned} Z_u &= \frac{4}{3\sqrt{\kappa_1}}\nabla u^{3/2} - \sqrt{\kappa_1 u}\,\nabla(K_s * w), \\ Z_w &= \frac{4}{3\sqrt{\kappa_2}}\nabla w^{3/2} - \sqrt{\kappa_2 w}\,\nabla(K_s * u). \end{aligned} \tag{4}$$

The energy dissipation and diffusion estimate are defined by

$$\begin{aligned} D(t) &= \int_{\mathbb{R}^3} (|Z_u|^2 + |Z_w|^2)(x,t)\,dx, \\ F(t) &= \frac{16}{9\kappa_1}\int_{\mathbb{R}^3} |\nabla u^{3/2}|^2\,dx + \frac{16}{9\kappa_2}\int_{\mathbb{R}^3} |\nabla w^{3/2}|^2\,dx. \end{aligned} \tag{5}$$

The gradients of the Sobolev powers and the weighted potential gradients vanish almost everywhere on the corresponding zero sets. These definitions do not presume an ordinary square-integrable density gradient.

**Theorem 2.1.** Under the stated assumptions, (1) admits a nonnegative global weak solution. For every finite $T > 0$, it has the regularity

$$u, w \in L^\infty(Q_T) \cap C([0,T]; L^1 \cap L^2), \qquad u^{3/2}, w^{3/2}, u^2, w^2 \in L^2(0,T; H^1). \tag{6}$$

For every $t \in [0,T]$, $\int_{\mathbb{R}^3} u(x,t)\,dx = M_u$ and $\int_{\mathbb{R}^3} w(x,t)\,dx = M_w$. The initial data are attained strongly in both spaces. Moreover,

$$\begin{aligned} &\operatorname{ess\,sup}_{0<t<T}(\|u(t)\|_\infty + \|w(t)\|_\infty) + \sup_{0\le t\le T}\left[\|u(t)\|_2^2 + \|w(t)\|_2^2 + \int_{\mathbb{R}^3} |x|^2 (u+w)(x,t)\,dx\right] \\ &\qquad + \int_0^T (D(t) + F(t))\,dt \le C_T, \end{aligned} \tag{7}$$

$$E(u(t), w(t)) + \int_a^t D(\sigma)\,d\sigma \le E(u(a), w(a)), \qquad 0 \le a \le t \le T. \tag{8}$$

The energy is continuous. The constant depends on $T, s, \kappa_1, \kappa_2$, initial masses, $L^\infty$ norms and second moments, and analytic inequality constants; it need not remain bounded as $T \to \infty$. The proof is given in Section 4.

### 3. Auxiliary lemmas and a priori estimates

Take radial $\rho \in C_c^\infty(B_1)$ with $\rho \ge 0$ and $\int_{\mathbb{R}^3} \rho dx = 1$. Set $\rho_\varepsilon(x) = \varepsilon^{-3}\rho(x/\varepsilon)$, $K_\varepsilon = \rho_\varepsilon * K_s$ and $a_\varepsilon(r) = r^2 + 2\varepsilon r$. Choose radial nonincreasing $\zeta \in C_c^\infty(B_2)$, equal to one on $B_1$, with $0 \le \zeta \le 1$, and put $\chi_R(x) = \zeta(x/R)$. For $f_0 = u_0, w_0$, choose

$$f_{0,\varepsilon} = \rho_\varepsilon * (\chi_{R_0+\varepsilon^{-1}} f_0), \qquad \| f_{0,\varepsilon} \|_1 \le \| f_0 \|_1, \quad \| f_{0,\varepsilon} \|_\infty \le \| f_0 \|_\infty. \tag{9}$$

Here $R_0$ is fixed so both truncated masses are positive. The smooth compactly supported data converge in every finite $L^p$ and in second moment. Their masses $M_{u,\varepsilon}, M_{w,\varepsilon}$ converge to $M_u, M_w$. The data do not depend on $R$.

**Lemma 3.1.** For each fixed $0 < \varepsilon \le 1$, the regularized Cauchy problem without spatial localization admits nonnegative global solutions $(u_\varepsilon, w_\varepsilon)$ with conserved masses and finite second moments on every finite time interval.

**Proof.** Localize by $\chi_R$ and extend $a_\varepsilon'$ smoothly and positively to the real line. Apply Matioc–Walker, Theorem 1.1.[14] In their notation choose $E_0 = (L^2)^2$, $E_1 = (H^2)^2$, $E_\vartheta = (H^{2\vartheta})^2$, and $(\gamma, \beta, \alpha, \xi, q) = (1/40, 4/5, 17/20, 9/10, 2)$. The principal operator generates an analytic semigroup on $L^2$ with domain $H^2$: its weighted realization is self-adjoint and its resolvent equation is coercive. Sobolev composition gives local Lipschitz dependence. The product $H^{4/5}H^{4/5} \subset H^{1/20}$ controls the quadratic gradient term, while $\nabla K_\varepsilon * : L^2 \to W^{2,\infty}$ controls the interaction. The required exponent condition is $(\xi - \alpha)q = 1/10 < 1 + \gamma - \alpha = 7/40$.

Negative-part and cutoff tests give nonnegativity and conserved approximate masses: the diffusion error is $O(L^{-1/2})$ by $L^2$ control, and the interaction error vanishes for $L > 2R$. Mass bounds the fields and their divergences independently of $R$; power testing gives an exponential $L^\infty$ bound. For either density $f$ with field $b$, put $q = a_\varepsilon(f)$. Testing by $-\Delta q$ with $a_\varepsilon'(f) \ge 2\varepsilon$ gives $f \in L_t^\infty H^1 \cap L_t^2 H^2$ and $f_t \in L_t^2 L^2$ on finite existence intervals, using $f = \sqrt{q + \varepsilon^2} - \varepsilon$ and $\| \nabla q \|_4^2 \le C \| q \|_\infty \| D^2 q \|_2$.

Guerand's Theorem 2.12 applies to the equation with bounded source $-\kappa f divb$ and uniformly elliptic coefficient $A = 2(f + \varepsilon)$, giving local $C^\nu$ regularity, $\nu > 0$.[15] Spatial Lipschitz continuity of $b$ gives $A$ and $-\kappa fb$ mean oscillation $O(r^\nu + r)$, hence Dini mean oscillation. Theorem 3.1 and estimate (3.17) of Dong et al. give local gradient bounds.[16] On fixed backward cylinders away from zero, constants are uniform in centre and terminal time. The global $H^1$ estimate controls the local gradient integral; taking the supremum over centres bounds the gradient globally.

For $v_h = [f(\cdot + he_k) - f]/h$, $k \in \{1,2,3\}$, the difference equation tested by $-\Delta v_h$ now yields

$$\frac{d}{dt} \| \nabla v_h \|_2^2 + 2\varepsilon \| \Delta v_h \|_2^2 \le C \| \nabla v_h \|_2^2 + C \| \Delta f \|_2^2 + C. \tag{10}$$

Choose a starting time in $(\tau/2, \tau)$ from the $L_t^2 H^2$ bound. Gronwall and $h \to 0$ give a uniform $H^2$ bound on $[\tau, T_*)$, where $T_* = \min(T, T_{\max})$. Constants depend on $\tau, \varepsilon, T$, not $R$ or $T_* - t$. Interpolation with $f_t \in L_t^2 L^2$ gives a $3/40$ time Hölder modulus in $H^{17/10}$. The uniform-continuity continuation criterion in Theorem 1.1(v)(a), with the whole interpolation space as state domain, excludes finite $T_{\max}$.[14]

For moments use $\varphi_L(x)=|x|^2/(1+|x|^2/L^2)$ with an outer compact cutoff. The bounds $|\nabla\varphi_L|\le 2\sqrt{\varphi_L}$ and $|\Delta\varphi_L|\le 14$ give a Gronwall estimate independent of $R,L$. Remove the outer cutoff first, then send $L\to\infty$. Local compactness and mass tails remove $\chi_R$. These fixed-$\varepsilon$ constants justify tests only. □

Set $I_\varepsilon=\int_{\mathbb{R}^3}u_\varepsilon(K_\varepsilon * w_\varepsilon)\,dx$ and

$$\begin{aligned} E_\varepsilon &= \frac{1}{\kappa_1}\int_{\mathbb{R}^3}[u_\varepsilon^2+2\varepsilon(u_\varepsilon\log u_\varepsilon-u_\varepsilon)]\,dx \\ &\quad+\frac{1}{\kappa_2}\int_{\mathbb{R}^3}[w_\varepsilon^2+2\varepsilon(w_\varepsilon\log w_\varepsilon-w_\varepsilon)]\,dx-I_\varepsilon, \\ J_\varepsilon &= \frac{1}{\kappa_1}\int_{\mathbb{R}^3}|x|^2\,u_\varepsilon\,dx+\frac{1}{\kappa_2}\int_{\mathbb{R}^3}|x|^2\,w_\varepsilon\,dx. \end{aligned} \tag{11}$$

Below $(v,\kappa,b)$ denotes either $(u_\varepsilon,\kappa_1,\nabla(K_\varepsilon * w_\varepsilon))$ or $(w_\varepsilon,\kappa_2,\nabla(K_\varepsilon * u_\varepsilon))$. Define $F_\varepsilon$ and $D_\varepsilon$ by adding the corresponding integrals of $|\nabla a_\varepsilon(v)/\sqrt{\kappa v}|^2$ and $|\nabla a_\varepsilon(v)/\sqrt{\kappa v}-\sqrt{\kappa v}\,b|^2$ over $\{x\in\mathbb{R}^3: v(x,t)>0\}$ with measure $dx$. We use $0\log 0=0$ and set the weighted integrands to zero on $\{v=0\}$.

**Lemma 3.2.** The solutions of Lemma 3.1 satisfy the energy identity (12) and estimate (16), with a constant independent of $0<\varepsilon\le 1$. In particular,

$$E_\varepsilon(t)+\int_a^t D_\varepsilon(\sigma)\,d\sigma=E_\varepsilon(a),\qquad 0\le a\le t. \tag{12}$$

**Proof.** Regularize entropy with $h_z''(r)=a_\varepsilon'(r)/[\kappa(r+z)]$, $0<z<\varepsilon$. Test by $h_z'(v)-h_z'(0)$; mass conservation removes the subtracted linear term. As $z\downarrow 0$, gradients are identified on the positive set and vanish on the zero set; moments control entropy tails. Truncate $K_\varepsilon$ by an even cutoff: symmetry pairs the two time derivatives, while the truncated gradients converge in $L^2$ and the endpoint energies converge by boundedness of $K_\varepsilon$ and mass. Combining gives (12), including the initial endpoint.

Put $U_\varepsilon=\|u_\varepsilon\|_2^2$, $V_\varepsilon=\|w_\varepsilon\|_2^2$ and $\beta=1-2s/3<1$. HLS and mass interpolation imply $I_\varepsilon\le C(U_\varepsilon V_\varepsilon)^{\beta/2}\le\eta(U_\varepsilon+V_\varepsilon)+C_\eta$. Negative entropy is controlled by mass and the already finite moment. Symmetrization yields

$$J_\varepsilon'=\frac{6}{\kappa_1}(U_\varepsilon+2\varepsilon M_{u,\varepsilon})+\frac{6}{\kappa_2}(V_\varepsilon+2\varepsilon M_{w,\varepsilon})-2(3-2s)I_\varepsilon+2\mathcal{R}_\varepsilon, \tag{13}$$

where the mollification remainder satisfies

$$\begin{aligned} \mathcal{R}_\varepsilon &= \int_{\mathbb{R}^3}\int_{\mathbb{R}^3}\int_{\mathbb{R}^3}u_\varepsilon(x)w_\varepsilon(y)\rho_\varepsilon(z)z\cdot\nabla K_s(x-y-z)\,dz\,dy\,dx, \\ |\mathcal{R}_\varepsilon| &\le \eta(U_\varepsilon+V_\varepsilon)+C_\eta\varepsilon^{3/(2s-1)},\quad 0<\varepsilon\le 1. \end{aligned} \tag{14}$$

Constants include masses and the mollifier's first moment. Gradient-kernel HLS with exponent $3/(1+s)$ gives a power $2(2-s)/3<1$; Young's inequality gives (14). With $\eta$ small, the entropy lower bound yields $\mathcal{G}_\varepsilon:=E_\varepsilon+3J_\varepsilon+C_0\ge c(U_\varepsilon+V_\varepsilon)+J_\varepsilon$. Combining (12) and (13), discarding $-6(3-2s)I_\varepsilon\le 0$, gives $\mathcal{G}_\varepsilon'+D_\varepsilon\le C\mathcal{G}_\varepsilon+C$. The initial values are uniformly bounded by (9) and moment convergence. Gronwall controls moments and quadratic norms jointly, without a prior uniform moment bound.

To recover the separate flux estimate, use $F_\varepsilon \le 2D_\varepsilon + 2W_\varepsilon$, where

$$W_\varepsilon = \kappa_1 \int_{\mathbb{R}^3} u_\varepsilon |\nabla(K_\varepsilon * w_\varepsilon)|^2 dx + \kappa_2 \int_{\mathbb{R}^3} w_\varepsilon |\nabla(K_\varepsilon * u_\varepsilon)|^2 dx. \tag{15}$$

Gradient-kernel HLS uses exponent $9/(1+4s)$. If this exceeds two, interpolation of the Sobolev power gives a total dissipation exponent $(7-8s)/7<1$; otherwise the mass and a priori $L^2$ bounds suffice. In both cases $W_\varepsilon \le F_\varepsilon/4 + C_T$. Consequently,

$$\sup_{0\le t\le T}(U_\varepsilon + V_\varepsilon + J_\varepsilon)(t) + \int_0^T (D_\varepsilon(t) + F_\varepsilon(t))\, dt \le C_T. \tag{16}$$

This proves the assertion. □

**Lemma 3.3.** For each finite $T>0$, both approximate densities are bounded in $L^\infty(Q_T)$ uniformly for $0<\varepsilon\le 1$.

**Proof.** For either triple defined above, compactly supported power tests followed by removal of the spatial cutoff give

$$\frac{1}{p}\frac{d}{dt}\| v\|_p^p + \frac{4(p-1)}{(p+1)^2}\| \nabla v^{(p+1)/2}\|_2^2 \le \frac{\kappa^2(p-1)}{4}\int_{\mathbb{R}^3} v^{p-1}|b|^2 dx, \qquad p\ge 2. \tag{17}$$

Only the test justification uses fixed-$\varepsilon$ $L^\infty$ and flux bounds; the estimates retain no such constants. Put $\alpha = 2s-1$, $z = 12/(3+2\alpha)$ and $\theta = (3-2\alpha)/5<1$. For the cubic test, HLS and Gagliardo–Nirenberg give $\| v\|_z \le C\| \nabla(v^2)\|_2^{\theta/2}\| v\|_2^{1-\theta}$. Adding the two cubic inequalities bounds their right-hand side by

$$C_T\| \nabla(u_\varepsilon^2)\|_2^\theta \| \nabla(w_\varepsilon^2)\|_2^\theta \le \frac{1}{4}\left(\| \nabla(u_\varepsilon^2)\|_2^2 + \| \nabla(w_\varepsilon^2)\|_2^2\right) + C_T. \tag{18}$$

This gives a uniform $L^3$ bound. Set $A_m = \sup_{0<\varepsilon\le 1}\max(\| u_\varepsilon\|_{L^\infty(0,T;L^m)}, \| w_\varepsilon\|_{L^\infty(0,T;L^m)})$. For $3\le m<3/\alpha$, take $p=2m-1$, $1/q = 1/m - \alpha/3$ and $r=(p-1)/(1-2/q)$. HLS bounds the right-hand integral in (17) by $C_{s,m}A_m^2\| v\|_r^{p-1}$. Interpolating $v^m$ between $L^1$ and $L^6$, with $m/r = \vartheta/6 + 1 - \vartheta$, gives $\| v\|_r^{p-1} \le C\| \nabla v^m\|_2^{2\theta} A_m^{me}$, where $e = 2(m-1)/m - 2\theta$ and $\theta = \vartheta(p-1)/(2m) = (3-2\alpha)/5 < 1$. Young absorption and $(2+me)/(1-\theta) = 2m$ yield

$$A_p^p \le M_0 A_0^{p-1} + C_{s,\kappa_1,\kappa_2,T,m}\, p^{2/(1-\theta)}\max(1, A_m^{2m}). \tag{19}$$

Here $M_0, A_0 \ge 1$ bound the initial masses and $L^\infty$ norms. Finitely many steps reach $m_* > 3/\alpha$; if an intermediate exponent hits the HLS endpoint, a slightly smaller target is used before the next step. Here $m_*' = m_*/(m_*-1)$. The near gradient kernel belongs to $L^{m_*'}$, and its far part is bounded. Mass and the $m_*$ bound yield an $\varepsilon$-independent bound for the interaction fields.

With the fields bounded, interpolate $v^m$ between $L^1$ and $L^6$ in (17). The absorption exponent is $a_m = 3(m-2)/(5m) < 3/5$, and the residual mass exponent is $b_m = (4m+2)/(5m)$. For $p=2m-1$, this gives

$$A_p^p \le M_0 A_0^{p-1} + C_T p^5 \max\left(1, A_m^{m(2m+1)/(m+3)}\right). \tag{20}$$

Here $C_T$ is independent of $m \ge m_*$: interpolation constants are bounded and $1/(1-a_m) \le 5/2$, so all factors are bounded by a fixed multiple of $p^5$. After the $p$th root, the old norm has exponent $\gamma_m = m(2m+1)/[(m+3)(2m-1)] < 1$, since $1-\gamma_m = (4m-3)/[(m+3)(2m-1)] > 0$. Along $m_n = 1+2^n(m_*-1)$, the logarithmic increments sum by $\sum_{n=0}^{\infty}(1+\log m_n)/m_n < \infty$. Thus both densities are uniformly bounded in $L^\infty(Q_T)$, including near zero. □

## 4. Global existence and energy inequality

**Proof.** Take the approximations from Lemma 3.1. Lemmas 3.2 and 3.3 provide uniform estimates on every finite time interval. Estimate (16) controls $\nabla u_\varepsilon^{3/2}$ and $\nabla w_\varepsilon^{3/2}$ in $L^2(Q_T)$. Cauchy–Schwarz with mass controls both fluxes in $L^2(0,T;L^1)$, so the time derivatives are bounded in $L^2(0,T;(W_0^{1,\infty})^*)$. Spatial translations of the Sobolev powers control density translations through $|a-b|^{3/2} \le |a^{3/2}-b^{3/2}|$ for $a,b \ge 0$. Spatial mollification and the time-derivative estimate give local compactness; second-moment tail estimates upgrade it to global strong $L^1(Q_T)$. Interpolation yields strong $L^p(Q_T)$ convergence for $1 \le p < 3$ and strong $L^1(Q_T)$ convergence of both diffusion squares.

Set $G_s = \nabla K_s$. Its near part belongs to $L^1$ and its far part to $L^2$. Hence

$$\| G_s * h\|_2 \le C(\| h\|_2 + \| h\|_1), \qquad h \in L^1(\mathbb{R}^3) \cap L^2(\mathbb{R}^3). \tag{21}$$

Apply this to $h = \rho_\varepsilon * w_\varepsilon - w$ and its counterpart. Strong $L^1(Q_T)$ convergence and the uniform mass bound give $\| h\|^2_{L^2(0,T;L^1)} \le C\| h\|_{L^1(Q_T)} \to 0$; strong $L^2(Q_T)$ convergence handles the other term. Thus the potential gradients converge in $L^2(Q_T)$, and their products with the densities converge in $L^1(Q_T)$. This identifies both equations in (2). Scalar time pairings recover the initial weak trace; expanding cutoffs and second-moment tail estimates preserve each mass. A diagonal extraction over integer time horizons gives one global pair.

For the first species, the complete vector in the approximate dissipation is

$$Z_{u,\varepsilon} = \frac{4}{3\sqrt{\kappa_1}}\nabla u_\varepsilon^{3/2} + \frac{4\varepsilon}{\sqrt{\kappa_1}}\nabla\sqrt{u_\varepsilon} - \sqrt{\kappa_1 u_\varepsilon}\nabla(K_\varepsilon * w_\varepsilon). \tag{22}$$

The middle term is bounded in $L^2(Q_T)$ and tends to zero distributionally. The first term converges weakly as the gradient of the density power, and the uniform density and interaction-field bounds identify the weighted potential gradient by strong convergence. The analogous statement holds for the second species. Lower semicontinuity therefore applies to the whole squared vector, not to individual signed cross terms. HLS gives convergence of interaction energies at almost every time; the additional entropy term vanishes by the uniform quadratic and moment bounds. Thus (8) initially holds on a common full-measure set of endpoints, also allowing the lower endpoint zero.

To include all endpoints, boundedness gives $u^2, w^2 \in L^2(0,T;H^1)$ and square-integrable fluxes. Hence $u_t, w_t \in L^2(0,T;H^{-1})$. Their continuous $H^{-1}$ representatives inherit the $L^1$, $L^\infty$ and moment bounds at every time by weak limits from regular times. Let $\mathcal{H} = H^{-1}$ with the Hilbert structure induced by $A = 1-\Delta$, and let $\Phi(v) = \frac{1}{3}\int_{\mathbb{R}^3} |v(x)|^3\, dx$ on $\mathcal{H} \cap L^3$, extended by infinity otherwise. Convexity and density show that $A(u^2) \in \partial_{\mathcal{H}}\Phi(u)$ almost everywhere. Both this subgradient and $u_t$ belong to $L^2(0,T;\mathcal{H})$. The subgradient inequalities and time difference quotients give the convex chain rule. At any endpoint $t$, average the inequality with competitor $u(t)$ over a one-sided interval; its error is bounded by the product of the $L^2$ norms of the subgradient and $u_t$ on that interval, hence tends to zero. This identifies the endpoint values and gives

$$\frac{1}{3}\int_{\mathbb{R}^3} u(x,t)^3\,dx - \frac{1}{3}\int_{\mathbb{R}^3} u(x,a)^3\,dx \\ = -\int_a^t\int_{\mathbb{R}^3}\Big([\nabla(u^2) - \kappa_1 u\nabla(K_s * w)]\cdot\nabla(u^2)\Big)(x,\sigma)\,dx\,d\sigma. \tag{23}$$

The same identity holds for $w$, with coefficient $\kappa_2$ and potential $K_s * u$. Continuity of the cubic norm and weak continuity imply strong $L^3$ continuity. The uniform moment tail then gives global $L^1$ continuity, and interpolation gives $L^2$ continuity, including the original initial data. Energy continuity follows from HLS; the dissipation primitive is continuous. Approximation of arbitrary endpoints by the common full-measure set proves (8) for every $0 \le a \le t \le T$. The bounds follow from Lemmas 3.2 and 3.3 by lower semicontinuity. This completes the proof of Theorem 2.1. □

## 5. Conclusion

Quadratic diffusion yields mass-conserving global weak solutions with finite-time boundedness and the energy inequality at all endpoints. Unequal orders, general diffusion powers, uniqueness and uniform-in-time bounds remain open here.